\documentclass[12pt]{article}

\usepackage{a4wide}
\usepackage{enumerate}
\usepackage[greek, frenchb, english]{babel}
\usepackage[T1]{fontenc}        
\usepackage{mathrsfs} \usepackage{amssymb}
\usepackage{amsmath}
\usepackage{amsfonts}
\usepackage[dvips
]{graphicx}
\usepackage{color} 
\usepackage[colorlinks=true, urlcolor=blue,linkcolor=blue, citecolor=blue]{hyperref}
\usepackage{aeguill}
\usepackage{rotating}
\usepackage{units}
\usepackage{graphicx}
\usepackage{hyperref}

\newtheorem{theo}{Theorem}[section]
\newtheorem{df}[theo]{Definition}
\newtheorem{lemme}[theo]{Lemma}

\newtheorem{prop}[theo]{Proposition}
\newtheorem{hyp}{Hypothesis}

\newcommand{\CQFD}{\hfill $\square$}

\def\real{\Bbb{R}}

\def\ee{\Bbb{E}}

\def\disp{\displaystyle}

\def\Tr{\mbox{Tr}}

\title{Rescaled weighted determinantal random balls}

\author{Adrien Clarenne\footnote{Univ Rennes, CNRS, IRMAR - UMR 6625, F-35000 Rennes, France. Email: adrien.clarenne@univ-rennes1.fr}}

\begin{document}

\maketitle

\begin{abstract}
We consider a collection of weighted Euclidian random balls in $\real^d$ distributed according a determinantal point process. We perform a zoom-out procedure by shrinking the radii while increasing the number of balls. We observe that the repulsion between the balls is erased and three different regimes are obtained, the same as in the weighted Poissonian case.
\end{abstract}


\section*{Introduction}

In this work, we give a generalization of the existing results concerning the asymptotics study of random balls model. The first results are obtained in 2007 by Kaj, Leskela, Norros and Schmidt in \cite{KLNS2007}. In their model, the balls are generated by an homogeneous Poisson point process on $\real^d\times\real_+$ (see \cite{DVJ1} for a general reference on point processes). In 2009, Breton and Dombry generalize this model adding in \cite{BD2009} a mark $m$ on the balls of the previous model and they obtain limit theorems on the so-called rescaled weighted random balls model. In 2010, Bierm\'e, Estrade and Kaj obtain in \cite{BEK2010} results performing for the first time a zoom-in scaling. In 2014, Gobard in his paper \cite{Gobard2014} extends the results of \cite{BD2009} considering inhomogeneous weighted random balls, and adding a dependence between the centers and the radii. The next step is to consider repulsion between the balls. In \cite{BCG2018}, Breton, Clarenne and Gobard give results on determinantal random balls model, but no weight are considered in their model. In this note, we consider weighted random balls generated by a non-stationary determinantal point process. To that purpose, we give an extension of the Laplace transform of determinantal processes allowing to compute Laplace transform with not necessarily compactly supported function, but with instead a condition of integrability. The main contributions of this note thus are a simplication of the proof of \cite{BCG2018} and the introduction of weights in the non-stationary determinantal random balls model.\\
This paper is organized as follows. In Section \ref{sec:model} and Section \ref{sec:asymptotics}, we give a description of the model and state our main results under the three different regimes. In Section \ref{sec:proof}, we give the Laplace transform of a determinantal point process for not compactly supported test functions and prove our results. Finally some technical results are gathered in the Appendix.

\section{Model}
\label{sec:model}
We consider a model of random balls in $\real^d$ constructed in the following way. The centers of the balls are generated by a determinantal point process (DPP) $\phi$ on $\real^d$ characterized by its kernel $K$ with respect to the Lebesgue measure. The motivation for considering such processes is that it introduces repulsion between the centers in agreement with various real model of balls. We assume that the map ${\bf K}$ given for all $f\in L^2(\real^d, dx)$ and $x\in \real^d$  by
\begin{equation}
\label{eq:KK}
{\bf K}f(x)=\int_{\real^d}K(x,y)f(y)\ dy
\end{equation}
satisfies the following hypothesis:
\begin{hyp} 
\label{hyp:det1}
The map ${\bf K}$ given in \eqref{eq:KK} is a bounded symmetric integral operator ${\bf K}$ from $L^2(\real^d,dx)$ into $L^2(\real^d,dx)$, with a continuous kernel $K$ with spectrum included in $[0,1[$. 
Moreover, ${\bf K}$ is locally trace-class, i.e. for all compact $\Lambda \subset E$, the restriction ${\bf K}_\Lambda$ of ${\bf K}$ on $L^2(\Lambda,\lambda)$ is of trace-class. 
\end{hyp}
Moreover, we also assume
\begin{eqnarray}
\label{eq:Kxx}
x \longmapsto K(x,x) \in L^{\infty}(\real^d).
\end{eqnarray}
These assumptions imply that $K(x,x)\geq0$.\\
At each center $x\in\real^d$, we attach two positive marks $r$ and $m$ independently. The first mark is interpreted as the radius and the second mark is the weight of the ball $B(x,r)$. The radii (resp. the weight) are independently and identically distributed according to $F$ (resp. according to $G$), assumed to admit a probability density $f$ (resp. a probability density $g$). We have a new point process $\Phi$ on $\real^d\times\real_+\times\real_+$ and according to Proposition A.7 in \cite{BCG2018}, $\Phi$ is a DPP on $\real^d\times \real_+\times \real_+$ with kernel
$$
\widehat{K}\big((x,r,m),(y,s,m')\big)= \sqrt{g(m)}\sqrt{f(r)}K(x,y)\sqrt{f(s)}\sqrt{g(m')},
$$
with respect to the Lebesgue measure.\\
\noindent Moreover, we suppose that the probability measure $G$ belongs to the normal domain of attraction of the $\alpha$-stable distribution $S_\alpha (\sigma,b,\tau)$ with $\alpha \in (1,2]$. Because $\alpha >1$, we can note that
\begin{equation}
\int_{\real_+} mG(dm)=\int_{\real_+}mg(m)dm<+\infty.
\end{equation}

\noindent In the sequel, we shall use the notation $\Phi$ both for the marked DPP
(i.e. the random locally finite collection of points  $(X_i, R_i, M_i)$) 
and for the associated random measure $\sum_{(X,R,M)\in \Phi}\delta_{(X,R,M)}$. 
We consider the contribution of the model in any suitable measure $\mu$ on $\real^d$ given by the following measure-indexed random field:  
\begin{equation}
\label{eq:procM}
M(\mu)=\int_{\real^d\times\real_+\times \real_+}m\mu\big(B(x,r)\big)\ \Phi(dx,dr,dm).
\end{equation}
However, in order to ensure that $M(\mu)$ in \eqref{eq:procM} is well defined, we restrain to measures $\mu$ with finite total variation  (see below Proposition~\ref{prop:Mesp}). 
In the sequel, ${\cal Z}(\real^d)$ stands the set of signed (Borelian) measures $\mu$ on $\real^d$ with finite total variation $\|\mu\|_{var}(\real^d)<+\infty$.
Moreover as in \cite{KLNS2007}, we assume the following assumption on the radius behaviour, for $d<\beta<2d$, 
\begin{equation}
\label{eq:tails}
f(r)\underset{r\to +\infty}{\sim} \frac{C_\beta}{r^{\beta+1}}, \qquad
r^{\beta+1} f(r)\leq C_0.
\end{equation}
Since $\beta>d$, condition \eqref{eq:tails} implies that the mean volume of the random ball is finite:
\begin{equation}
\label{eq:cfvol}
v_d\int_0^{+\infty}r^df(r)\ dr<+\infty,
\end{equation}
where $v_d$ is the Lebesgue measure of the unit ball of $\real^d$.
On the contrary, $\beta<2d$ implies that $F$ does not admit a moment of order $2d$ and the volume of the balls has an infinite variance. 
The asymptotics condition in \eqref{eq:tails} is of constant use in the following.  

\begin{prop}
\label{prop:Mesp}
Assume \eqref{eq:tails} is in force. 
For all $\mu \in {\cal Z}(\real^d)$, $\ee\left[ M(\vert\mu\vert)  \right]<+ \infty$. 
As a consequence, $M(\mu)$ in \eqref{eq:procM} is almost surely well defined for all $\mu\in{\cal Z}(\real^d)$.
\end{prop}
\begin{Proof}
The proof follows the same lines as that of Proposition 1.1 in \cite{BCG2018}, replacing $K(0)$ by $K(x,x)$ and controlling it thanks to Hypothesis \eqref{eq:Kxx}.
\CQFD
\end{Proof}

\section{Asymptotics and main results}
\label{sec:asymptotics}

The zooming-out procedure acts accordingly both on the centers and on the radii. First, a scaling $S_\rho: r\mapsto \rho r$ of rate $\rho\in(0,1]$ changes balls $B(x,r)$ into $B(x,\rho r)$; 
this scaling changes the distribution $F$ of the radius into $F_\rho=F\circ S_\rho^{-1}$.
Second, the intensity of the centers is simultaneously adapted; 
to do this, we introduce actually a family of new kernels $K_\rho$, $\rho\in]0,1]$, that we shall refer to as scaled kernels, and we denote by $\phi_\rho$ the DPP with kernel $K_\rho$ (with respect to the Lebesgue measure).The zoom-out procedure consists now in introducing the family of DPPs $\phi_\rho$, $\rho\in]0,1]$, with kernels $K_\rho$ with respect to the Lebesgue measure satisfying
\begin{equation}\label{chap2:hypkrho0}
K_\rho(x,x)\underset{\rho \rightarrow 0}{\sim}\lambda (\rho)K(x,x),
\end{equation}
with $\lim_{\rho\to 0} \lambda (\rho)= +\infty$. 
We also suppose
\begin{equation}
\label{eq:lambdarho2}
\sup_{x\in\real^d} K_\rho(x,x)\leq \lambda (\rho)\sup_{x\in\real^d} K(x,x),
\end{equation}
and observe that with \eqref{eq:Kxx} and \eqref{eq:lambdarho2}, Proposition A.6 in \cite{BCG2018} gives the following uniform bound
\begin{equation}
\label{eq:controlintKrho2} 
\sup_{x\in\real^d}\int_{\real^d}\big|K_\rho(x,y) \big|^2 \ dy
\underset{\rho \rightarrow 0}{=}\mathcal{O}\big(\lambda(\rho)\big).
\end{equation}

\medskip\noindent
The zoom-out procedure consists in considering a new marked DPP $\Phi_\rho$ on $\real^d\times \real_+\times \real_+$ with kernel:
$$
\widehat{K}_\rho\big((x,r,m),(y,s,m')\big)
=\sqrt{g(m)}\sqrt{\frac{f(r/\rho)}{\rho}}K_\rho(x,y)\sqrt{\frac{f(s/\rho)}{\rho}}\sqrt{g(m')},
$$
with respect to the Lebesgue measure. 
The so-called scaled version of $M(\mu)$ is then the field
$$
M_{\rho}(\mu)=\int_{\real^d\times\real_+\times \real_+} m\mu\big(B(x,r)\big)\ \Phi_\rho(dx,dr,dm).
$$
In the sequel, we are interested in the fluctuations of $M_{\rho}(\mu)$ with respect to its expectation 
$$
\ee\big[M_{\rho}(\mu)\big]=\int_{\real^d\times\real_+\times \real_+}m\mu\big(B(x,r)\big)\ K_\rho(x,x)\frac{f(r/\rho)}{\rho}g(m) dxdrdm
$$
and we introduce
\begin{equation}
\label{eq:procMM}
\widetilde M_\rho(\mu)
=M_\rho(\mu)-\ee\big[M_\rho(\mu)\big]
=\int_{\real^d\times\real_+\times \real_+} m\mu\big(B(x,r)\big)\ \widetilde\Phi_\rho(dx,dr,dm),
\end{equation}
 where $\widetilde\Phi_\rho$ stands for the compensated random measure associated to $\Phi_\rho$.\\ 

We introduce a subspace $\mathcal{M}_{\alpha,\beta}\subset \mathcal{Z}$ on which we will investigate the convergence of the random field $M_\rho(\mu)$. The next definition comes from \cite{BD2009}. 

\begin{df}
For $1<\alpha\leq 2$ and $\beta >0$, we denote by $\mathcal{M}_{\alpha,\beta}$ the subset of measures $\mu\in\mathcal{Z}(\real^d)$ satisfying for some finite constant $C_\mu$ and some $0<p<\beta<q$ :
\[\int_{\real^d}\left|\mu(B(x,r))\right|^\alpha dx \leq C_\mu (r^p \wedge r^q)\]
where $a\wedge b=\min(a,b)$.
\end{df}
We denote by $\mathcal{M}_{\alpha,\beta}^+$ the space of positive measures $\mu\in\mathcal{M}_{\alpha,\beta}$.
Now, we can state the main result of this note. The proof consists in a combination of the arguments of \cite{BCG2018} and \cite{Gobard2014}. It is given in Section \ref{sec:proof} where for some required technical points, it is referred to \cite{BCG2018} and \cite{BD2009}.
\newpage
\begin{theo} \label{main theorem}
Assume \eqref{eq:tails} and $\phi_\rho$ is a DPP with kernel satisfying \eqref{eq:Kxx}, \eqref{chap2:hypkrho0}, \eqref{eq:lambdarho2} and Hypothesis~\ref{hyp:det1} for its associated operator ${\bf K}_\rho$ in \eqref{eq:KK}.

\begin{enumerate}[(i)]

\item{Large-balls scaling}: Assume $\lambda(\rho)\rho^\beta\to +\infty$ and set $n(\rho)=\big(\lambda(\rho)\rho^\beta\big)^{1/\alpha}$. 
Then, $\widetilde M_\rho(\cdot)/n(\rho)$ converges in the fdd sense on ${\cal M}_{\alpha,\beta}^+$ to $W_\alpha (\cdot)$  where
$$
W_\alpha(\mu)=\int_{\real^d\times\real_+} \mu\big(B(x,r)\big)\ M_\alpha(dx,dr)
$$
is a stable integral with respect to the $\alpha$-stable random measure $M_\alpha$ with control measure $\sigma^\alpha K(x,x)C_\beta r^{-\beta-1}\ dxdr$ and constant skewness function $b$ given in the domain of attraction of $G$.

\item{Intermediate scaling}: Assume $\lambda(\rho)\rho^\beta\to a^{d-\beta}\in]0,+\infty[$ and set $n(\rho)=1$. 
Then, $\widetilde M_\rho(\cdot)/n(\rho)$ converges in the fdd sense on ${\cal M}_{\alpha,\beta}^+$ to $\widetilde P\circ D_a$ where  
$$
\widetilde P(\mu)=\int_{\real^d\times\real_+\times\real_+} m\mu\big(B(x,r)\big)\ \widetilde \Pi(dx,dr,dm)
$$
with $\widetilde \Pi$ a (compensated) PPP with compensator measure $K(x,x)C_\beta r^{-\beta-1}\ dxdrG(dm)$ and $D_a$ is the dilatation defined by $(D_a\mu)(B)=\mu(a^{-1}B)$.

\item{Small-balls scaling}: Suppose $\lambda(\rho)\rho^\beta \to 0$ when $\rho \to 0$ for $d<\beta<\alpha d$ and set $n(\rho)=(\lambda(\rho)\rho^\beta)^{1/\gamma}$ with $\gamma=\beta/d \in ]1,\alpha[$. 
Then, the field $n(\rho)^{-1}\widetilde{M}_\rho(\cdot)$ converges when $\rho\to 0$ in the finite-dimensional distributions sense in $L^1(\real^d)\cap L^2(\real^d)\cap\left\lbrace \mu \geq 0\right\rbrace$ to $Z_\gamma(\cdot)$ where
$$
Z_\gamma(\mu)=\int_{\real^d} \phi(x)\ M_\gamma(dx) \quad \mbox{ for } \quad\mu(dx)=\phi(x)dx \mbox{ with } \phi\in L^1(\real^d)\cap L^2(\real^d),~\phi \geq 0,
$$
is a stable integral with respect to the $\gamma$-stable random measure $M_\gamma $ with control measure $\sigma_\gamma K(x,x) dx$ where
$$
\sigma_\gamma^\gamma=\frac{ C_\beta v_d^\gamma}d \int_{0}^{+\infty} \frac{1-\cos(r)}{r^{1+\gamma}}dr\int_0^{+\infty}m^\gamma G(dm),
$$
and constant unit skewness.
\end{enumerate}
\end{theo}
Here, and in the sequel, we follow the notations of the standard reference \cite{ST1994} for stable random variables and integrals.


\section{Proof}
\label{sec:proof}

To investigate the behaviour of $\widetilde{M}_\rho(\mu)$ in the determinantal case, we use the Laplace transform of determinantal measures. An explicit expression is well known when the test functions are compactly supported, see Theorem A.4 in \cite{BCG2018}. However, in our situation, the test functions $(x,r,m)\longmapsto m\mu(B(x,r))$ are not compactly supported on $\real^d\times\real_+\times\real_+$ for $\mu\in\mathcal{M}_{\alpha,\beta}^+$. In order to overpass this issue we use Proposition \ref{Laplace} below for the Laplace transform of determinantal measures with non-compactly supported test functions, but with a condition of integration with respect to the kernel of the determinantal process (see~\eqref{cond:integrability}). \\
In addition to generalizing the model studied in \cite{BCG2018} by adding a weight, the following proposition has the further consequence of simplifying the proofs of the results in \cite{BCG2018}, since there is no more need to study the truncated model and obtain uniform convergence to exchange the limit in $R$, the truncation parameter and the limit in $\rho$, the scaling parameter. 

\begin{prop}\label{Laplace}
Let $\Phi$ a determinantal point process on a locally compact Polish space E with a continuous kernel $K$ such that the associated operator $\textbf{K}$ satisfies Hypothesis~\ref{hyp:det1}. We also suppose that the kernel $K\left[1-e^{-h}\right]\in L^2(E\times E)$.  Let $h$ be a nonnegative function such that 
\begin{equation}\label{cond:integrability}
	\int_E \left(1-e^{-h(x)}\right) K(x,x)dx <+\infty .
\end{equation}
Then $\mathbf{K}\left[1-e^{-h}\right]$ is a trace-class operator with $$\mathrm{Tr} \left(\mathbf{K}\left[1-e^{-h}\right]\right)=\int_E \left(1-e^{-h(x)}\right) K(x,x)dx$$ and we have
\begin{equation}\label{laplace eq}
	\mathbb{E}\left[\exp\left(-\int_E h(x)\Phi(dx)\right)\right]=\exp\left(-\sum_{n= 1}^{\infty} \frac{1}{n}\mathrm{Tr}\Big(\mathbf{K}\left[1-e^{-h}\right]^n\Big)\right),
\end{equation}
where $\mathbf{K}\left[1-e^{-h}\right]$ is the operator with kernel 
$$K\left[1-e^{-h}\right](x,y)=\sqrt{1-e^{-h(x)}}K(x,y) \sqrt{1-e^{-h(y)}}.$$
\end{prop}

\noindent\begin{Proof}
Expression (\ref{laplace eq}) is known to be true when $h$ has a compact support (see Theorem A.4 and equation (37) in \cite{BCG2018}), but it is not the case here. Let $(h_p)_{p\in\mathbb{N}}$ a non-decreasing sequel of positive functions with compact support defined by
\[h_p(x)=h(x)\textbf{1}_{B(0,p)}(x).\]	
Thanks to Theorem A.4 in \cite{BCG2018}, we have for all $p\in\mathbb{N}$:
\begin{equation}\label{laplacelimit}
	\mathbb{E}\left[\exp\left(-\int_E h_p(x)\Phi(dx)\right)\right]=\exp\left(-\sum_{n= 1}^{\infty} \frac{1}{n}\Tr\Big(\mathbf{K}\left[1-e^{-h_p}\right]^n\Big)\right).
\end{equation}
\underline{First step}: To begin, we prove that \[\displaystyle\underset{p\rightarrow +\infty}{\lim}\mathbb{E}\left[\exp\left(-\int_E h_p(x)\Phi(dx)\right)\right]=\mathbb{E}\left[\exp\left(-\int_E h(x)\Phi(dx)\right)\right].\] \\
If we denote by $M_p(E)$ the space of all point measures defined on $E$, we have 
$$\mathbb{E}\left[\exp\left(-\int_E h_p(x)\Phi(dx)\right)\right]=\int_{M_p(E)}\exp\left(-\int_E h_p(x)m(dx)\right)\mathbb{P}_\Phi (dm).$$
Because $p\longmapsto h_p$ is increasing, by monotone convergence we have 
\[\underset{p\rightarrow +\infty}{\lim}\int_E h_p(x)m(dx) = \int_E h(x)m(dx).\]
To finish we apply the dominated convergence theorem:
\begin{itemize}
\item $\displaystyle\underset{p\rightarrow +\infty}{\lim}\exp\left(-\int_E h_p(x)m(dx)\right)=\exp\left(-\int_E h(x)m(dx)\right)$,
\item $\displaystyle\left|\exp\left(-\int_E h_p(x)m(dx)\right)\right|\leq 1$ because $h_p \geq 0$, and 1 is integrable on $M_p(E)$ with respect to $\mathbb{P}_\Phi$.\\
\end{itemize}
\underline{Second step}: We prove that $\mathbf{K}\left[1-e^{-h}\right]$ is a trace-class operator.\\\\
Let $(e_n)_{n\in\mathbb{N}}$ be an orthonormal basis of $L^2(E)$.\\
$\mathbf{K}\left[1-e^{-h}\right]$ is trace-class if $\disp \sum_{n=0}^{+\infty} \left\langle \mathbf{K}\left[1-e^{-h}\right](e_n),e_n \right\rangle <+\infty$ and in this case we have $\disp\Tr \left(\mathbf{K}\left[1-e^{-h}\right]\right)=\sum_{n=0}^{+\infty} \left\langle \mathbf{K}\left[1-e^{-h}\right](e_n),e_n \right\rangle.$\\
We begin by proving that for all $n\in\mathbb{N}$, $$\underset{p\rightarrow +\infty}{\lim}\left\langle \mathbf{K}\left[1-e^{-h_p}\right](e_n),e_n \right\rangle=\left\langle \mathbf{K}\left[1-e^{-h}\right](e_n),e_n \right\rangle.$$
\begin{align*}
\left\langle \mathbf{K}\left[1-e^{-h_p}\right](e_n),e_n \right\rangle=\int_{E\times E}\sqrt{1-e^{-h_p(x)}} K(x,y) \sqrt{1-e^{-h_p(y)}}e_n(x)e_n(y)dxdy
\end{align*}
We apply the dominated convergence theorem:
\begin{itemize}
\item Computation of the limit when $p\rightarrow +\infty$:
\begin{align*}
\underset{p\rightarrow +\infty}{\lim}\sqrt{1-e^{-h_p(x)}} K(x,y) & \sqrt{1-e^{-h_p(y)}}e_n(x)e_n(y)\\
&=\sqrt{1-e^{-h(x)}} K(x,y) \sqrt{1-e^{-h(y)}}e_n(x)e_n(y)
\end{align*}
\item Domination:
\begin{align*}
&\left|\sqrt{1-e^{-h_p(x)}} K(x,y) \sqrt{1-e^{-h_p(y)}}e_n(x)e_n(y)\right|\\
&\hspace{4cm}\leq\sqrt{1-e^{-h(x)}} \left|K(x,y)\right| \sqrt{1-e^{-h(y)}}\left|e_n(x)\right|\left|e_n(y)\right|
\end{align*}
which is integrable on $E\times E$ because:
\begin{align*}
&\int_{E\times E}\sqrt{1-e^{-h(x)}} \left|K(x,y)\right| \sqrt{1-e^{-h(y)}}\left|e_n(x)\right|\left|e_n(y)\right|dxdy \\
&\leq \sqrt{\int_{E\times E}(1-e^{-h(x)}) \left|K(x,y)\right|^2 (1-e^{-h(y)})dxdy} \times \sqrt{\int_{E\times E}\left|e_n(x)\right|^2\left|e_n(y)\right|^2dxdy}\\
&\leq \int_{E}(1-e^{-h(x)}) K(x,x)dx  <+\infty
\end{align*}
because $\left|K(x,y)\right|^2\leq K(x,x)K(y,y)$ and $\disp\int_{E}\left|e_n(x)\right|^2dx=1$.
\end{itemize}
Thus we have:
\begin{align*}
\sum_{n=0}^{+\infty} \left\langle \mathbf{K}\left[1-e^{-h}\right](e_n),e_n \right\rangle&=\sum_{n=0}^{+\infty}\underset{p\rightarrow +\infty}{\mathrm{lim~inf}} \left\langle \mathbf{K}\left[1-e^{-h_p}\right](e_n),e_n \right\rangle\\
&\underset{(a)}{\leq}\underset{p\rightarrow +\infty}{\mathrm{lim~inf}}\sum_{n=0}^{+\infty} \left\langle \mathbf{K}\left[1-e^{-h_p}\right](e_n),e_n \right\rangle \\
&=\underset{p\rightarrow +\infty}{\mathrm{lim~inf}}~ \Tr \left(\mathbf{K}\left[1-e^{-h_p}\right]\right)
\end{align*}
where (a) is allowed by the Fatou lemma because $\left\langle \mathbf{K}\left[1-e^{-h_p}\right](e_n),e_n \right\rangle\geq 0$.\\
We have 
\begin{align*}
\Tr \left(\textbf{K}\left[1-e^{-h_p}\right]\right)=\int_{E} \left(1-e^{-h_p(x)}\right)K(x,x)dx.
\end{align*}
We apply the dominated convergence theorem: 
\begin{enumerate}
\item $\underset{p\rightarrow +\infty}{\lim}\left(1-e^{-h_p(x)}\right)K(x,x)=\left(1-e^{-h(x)}\right)K(x,x)$,
\item For all $p\geq 0$ we have:
\begin{equation*}
\left| \left(1-e^{-h_p(x)}\right)K(x,x) \right| \leq \left(1-e^{-h(x)}\right)K(x,x) 
\end{equation*}
which is integrable on $E$ by hypothesis \eqref{cond:integrability}.
\end{enumerate}
Thus
\begin{equation*}
\underset{p\rightarrow +\infty}{\lim}\int_E \left(1-e^{-h_p(x)}\right) K(x,x)dx=\int_E \left(1-e^{-h(x)}\right) K(x,x)dx.
\end{equation*}
So 
\[\sum_{n=0}^{+\infty} \left\langle \mathbf{K}\left[1-e^{-h}\right](e_n),e_n \right\rangle\leq \int_E \left(1-e^{-h(x)}\right) K(x,x)dx<+\infty\]
which proves that $\mathbf{K}\left[1-e^{-h}\right]$ is trace-class.\\

\noindent \underline{Third step}: We prove that $\disp \Tr \left(\mathbf{K}\left[1-e^{-h}\right]\right)=\int_E \left(1-e^{-h(x)}\right) K(x,x)dx.$\\
Because $\mathbf{K}\left[1-e^{-h}\right]$ is trace-class and $K\left[1-e^{-h}\right]\in L^2(E\times E)$, Lemma 4.2.2 in \cite{HKPV2009} ensures that 
\begin{equation}\label{chap2:decompnoyau}
\sqrt{1-e^{-h(x)}} K(x,y) \sqrt{1-e^{-h(y)}}=\sum_{i=1}^{+\infty}\lambda_i \varphi_i(x)\varphi_i(y) 
\end{equation}
where $(\varphi_i)$ is an orthonormal basis of $L^2(E)$. Thus, thanks to \eqref{chap2:decompnoyau} we have:
\begin{align}\label{chap2:premiereinterversion}
 \left\langle \mathbf{K}\left[1-e^{-h}\right](e_n),e_n \right\rangle&=\int_{E\times E}\sqrt{1-e^{-h(x)}}K(x,y) \sqrt{1-e^{-h(y)}}e_n(x)e_n(y)dxdy \nonumber\\
 &=\int_{E\times E}\sum_{i=1}^{+\infty}\lambda_i  \varphi_i(x)\varphi_i(y)e_n(x)e_n(y)dxdy
\end{align}
We can invert the sum and the integral if $$\sum_{i=1}^{+\infty} \int_{E\times E} \left| \lambda_i \varphi_i(x)\varphi_i(y)  e_n(x)e_n(y)\right|dxdy<+\infty.$$
\begin{align*}
\sum_{i=1}^{+\infty}\int_{E\times E} \left| \lambda_i \varphi_i(x)\varphi_i(y)  e_n(x)e_n(y)\right|dxdy&=\sum_{i=1}^{+\infty}\int_{E\times E}  \lambda_i  \left|\varphi_i(x)\right|\left|\varphi_i(y)\right|  \left|e_n(x)\right|\left|e_n(y)\right|dxdy \\
&=\sum_{i=1}^{+\infty}\lambda_i \left(\int_{E}\left|\varphi_i(x)\right|\left|e_n(x)\right|dx\right)^2\\
&\underset{(1)}{\leq} \sum_{i=1}^{+\infty}\lambda_i \left(\int_{E} \left|\varphi_i(x)\right|^2dx\right)\left(\int_E \left|e_n(x)\right|^2dx \right)\\
&\underset{(2)}{=} \sum_{i=1}^{+\infty}\lambda_i \int_{E} \left|\varphi_i(x)\right|^2dx\\
&\underset{(3)}{=} \int_{E} \left(\sum_{i=1}^{+\infty}\lambda_i  \left|\varphi_i(x)\right|^2\right)dx\\
&\underset{(4)}{=} \int_{E}(1-e^{-h(x)})K(x,x)dx <+\infty.
\end{align*}
The inequality (1) is the Cauchy-Schwarz inequality, the equality (2) is because 
$$\disp \int_E \left|e_n(x)\right|^2dx = \left\|e_n\right\|_2^2=1,$$
the inversion of the sum and the integral in (3) is allowed because all the terms are non-negative and the equality (4) stands because $K$ is continuous. \\
So we can invert the integral and the sum in \eqref{chap2:premiereinterversion} to obtain
\begin{align*}
 \left\langle \mathbf{K}\left[1-e^{-h}\right](e_n),e_n \right\rangle&=\sum_{i=1}^{+\infty}\int_{E\times E}\lambda_i  \varphi_i(x)\varphi_i(y)e_n(x)e_n(y)dxdy\\
 &=\sum_{i=1}^{+\infty}\lambda_i \left(\int_{E}  \varphi_i(x)e_n(x)dx\right)^2\\
 &=\sum_{i=1}^{+\infty}\lambda_i \left\langle  \varphi_i,e_n \right\rangle^2
\end{align*}
and the result is proved because
\begin{align*}
 \sum_{n=0}^{+\infty} \left\langle \mathbf{K}\left[1-e^{-h}\right](e_n),e_n \right\rangle &= \sum_{n=0}^{+\infty} \sum_{i=1}^{+\infty}\lambda_i \left\langle \varphi_i,e_n \right\rangle^2 \\
 &\underset{(a)}{=} \sum_{i=1}^{+\infty} \lambda_i \sum_{n=0}^{+\infty}\left\langle  \varphi_i,e_n \right\rangle^2 \\
 &=\sum_{i=1}^{+\infty} \lambda_i \left\| \varphi_i \right\|_2^2 \\
 &=\sum_{i=1}^{+\infty} \lambda_i \int_E \left|\varphi_i(x)\right|^2 dx\\
&\underset{(a)}{=} \int_E\sum_{i=1}^{+\infty} \lambda_i  \left|\varphi_i(x)\right|^2 dx\\
&=\int_{E}(1-e^{-h(x)})K(x,x)dx <+\infty,
\end{align*}
where the two equalities (a) are allowed because the terms are non-negative.\\
Then, $\mathbf{K}\left[1-e^{-h}\right]$ is a trace-class operator with \[\disp\Tr \left(\mathbf{K}\left[1-e^{-h}\right]\right)=\int_E \left(1-e^{-h(x)}\right) K(x,x)dx.\]
The computations above gives in particular
\[\underset{p\rightarrow +\infty}{\lim} \Tr \left(\textbf{K}\left[1-e^{-h_p}\right]\right)=\Tr \left(\textbf{K}\left[1-e^{-h}\right]\right).\]
\noindent 
\underline{Fourth step}: Now, we have to take the limit in the right term of equation (\ref{laplacelimit}). \\\\
First, we prove that we can exchange the limit and the infinite sum. To do that, we show that the sum normally converges. For all $n,p\geq 0$:
\begin{align*}
	\left|\Tr\Big(\textbf{K}\left[1-e^{-h_p}\right]^n\Big)\right|&\leq \left\| \textbf{K}\left[1-e^{-h_p}\right]\right\|^{n-1} \Tr\Big(\textbf{K}\left[1-e^{-h_p}\right]\Big)\\
	&\leq \left\| \textbf{K}\left[1-e^{-h}\right]\right\|^{n-1} \Tr\Big(\textbf{K}\left[1-e^{-h}\right]\Big),
\end{align*}
where the second inequality stands thanks to Lemma \ref{ineg norme} in Appendix \ref{sec:laplacelemme} and because $K(x,x)\geq 0$, 
\begin{align}
\Tr\Big(\textbf{K}\left[1-e^{- h_p}\right]\Big)&=\int_E \left(1-e^{-h_p(x)}\right) K(x,x)dx \nonumber\\
&\leq \int_E \left(1-e^{-h(x)}\right) K(x,x)dx \nonumber\\
&=\Tr\Big(\textbf{K}\left[1-e^{- h}\right]\Big),
\label{trace-control}
\end{align}
which is finite as proved above.\\
\noindent Let now prove that $\left\|\textbf{K}\left[1-e^{-h}\right]\right\|<1$.
\begin{align*}
\left\|\textbf{K}\left[1-e^{-h}\right]\right\|&=\underset{\|g\|_2=1}{\sup}\left\langle \textbf{K}\left[1-e^{-h}\right](g),g \right\rangle \\
&=\underset{\|g\|_2=1}{\sup}\hspace{0.2cm}\underset{p\rightarrow +\infty}{\lim}\left\langle \textbf{K}\left[1-e^{-h_p}\right](g),g \right\rangle \\
&\leq \underset{p\rightarrow +\infty}{\lim}\hspace{0.2cm}\underset{\|g\|_2=1}{\sup}\left\langle \textbf{K}\left[1-e^{-h_p}\right](g),g \right\rangle \\
&=\underset{p\rightarrow +\infty}{\lim}\hspace{0.2cm}\underset{\|g\|_2=1}{\sup}\left\langle \textbf{K}\left(\sqrt{1-e^{-h_p}}g\right),\sqrt{1-e^{-h_p}}g \right\rangle
\end{align*}
Because $\sqrt{1-e^{-h_p}}g$ has a compact support, 
\begin{equation*}
 \left\langle \textbf{K}\left(\sqrt{1-e^{-h_p}}g\right),\sqrt{1-e^{-h_p}}g \right\rangle=\left\langle \textbf{K}_{\left|B(0,p)\right.}\left(\sqrt{1-e^{-h_p}}g\right),\sqrt{1-e^{-h_p}}g \right\rangle,
\end{equation*}
where $\textbf{K}_{\left|B(0,p)\right.}$ is the restriction of $\textbf{K}$ on $L^2(B(0,p))$. Since $\textbf{K}_{\left|B(0,p)\right.}$ is trace-class, we have 
\begin{align*}
\left\langle \textbf{K}_{\left|B(0,p)\right.}\left(\sqrt{1-e^{-h_p}}g\right),\sqrt{1-e^{-h_p}}g \right\rangle &\leq \lambda_{p}^{\mathrm{max}} \left\|\sqrt{1-e^{-h_p}}g \right\|^2 \\
&\leq \lambda_{p}^{\mathrm{max}} \\
&\leq \lambda^{\mathrm{max}},
\end{align*} 
where $\lambda_{p}^{\mathrm{max}}$ (resp. $\lambda^{\mathrm{max}}$) is the greatest eigenvalue of $\textbf{K}_{\left|B(0,p)\right.}$ (resp $\textbf{K}$). Then
\[\underset{\|g\|_2=1}{\sup}\left\langle \textbf{K}\left(\sqrt{1-e^{-h_p}}g\right),\sqrt{1-e^{-h_p}}g \right\rangle\leq \lambda^{\mathrm{max}}\]
and 
\[\underset{p\rightarrow +\infty}{\lim}\underset{\|g\|_2=1}{\sup}\left\langle \textbf{K}\left(\sqrt{1-e^{-h_p}}g\right),\sqrt{1-e^{-h_p}}g \right\rangle\leq \lambda^{\mathrm{max}}.\]
Thus,
\[\left\|\textbf{K}\left[1-e^{-h}\right]\right\|\leq \lambda^{\mathrm{max}} <1\]
which proves the result.\\\\
So we have an upper bound of $\displaystyle\frac{1}{n}\Tr\Big( \textbf{K}\left[1-e^{-h_p}\right]^n\Big)$ independent of $p$ which is summable so we can exchange the limit and the sum.\\\\
Secondly, we prove that $\underset{p\rightarrow +\infty}{\lim}\Tr\Big(\textbf{K}\left[1-e^{-h_p}\right]^n\Big)=\Tr\Big(\textbf{K}\left[1-e^{-h}\right]^n\Big)$.\\
For all $n\geq 0$:
\begin{align*}
	&\left|\Tr\Big( \textbf{K}\left[1-e^{-h}\right]^n\Big)-\Tr\Big(\textbf{K}\left[1-e^{-h_p}\right]^n\Big)\right|=\left| \Tr\Big( \textbf{K}\left[1-e^{-h}\right]^n- \textbf{K}\left[1-e^{-h_p}\right]^n\Big) \right| \\
	&\hspace*{2.2cm}= \left| \Tr\left( \left(\textbf{K}\left[1-e^{-h}\right]-\textbf{K}\left[1-e^{-h_p}\right]\right)\sum_{k=0}^{n-1} \textbf{K}\left[1-e^{-h}\right]^k \textbf{K}\left[1-e^{-h_p}\right]^{n-1-k}\right) \right| \\
	&\hspace*{2.2cm}\leq \sum_{k=0}^{n-1}\Tr \left(\left(\textbf{K}\left[1-e^{-h}\right]-\textbf{K}\left[1-e^{-h_p}\right]\right)\textbf{K}\left[1-e^{-h}\right]^{k} \textbf{K}\left[1-e^{-h_p}\right]^{n-1-k} \right). 
\end{align*}
Moreover we have the following inequalities
\begin{align*}
\Tr &\left(\left(\textbf{K}\left[1-e^{-h}\right]-\textbf{K}\left[1-e^{-h_p}\right]\right)\textbf{K}\left[1-e^{-h}\right]^{k} \textbf{K}\left[1-e^{-h_p}\right]^{n-1-k} \right)\\
&\hspace*{2.2cm}\leq  \Tr \left(\textbf{K}\left[1-e^{-h}\right]-\textbf{K}\left[1-e^{-h_p}\right]\right)\left\| \textbf{K}\left[1-e^{-h}\right]^{k} \textbf{K}\left[1-e^{-h_p}\right]^{n-1-k}\right\| \\
&\hspace*{2.2cm}\leq  \Tr \left(\textbf{K}\left[1-e^{-h}\right]-\textbf{K}\left[1-e^{-h_p}\right]\right)\left\| \textbf{K}\left[1-e^{-h}\right]\right\|^{k} \left\|\textbf{K}\left[1-e^{-h_p}\right]\right\|^{n-1-k}\\
&\hspace*{2.2cm}\leq\Tr \left(\textbf{K}\left[1-e^{-h}\right]-\textbf{K}\left[1-e^{-h_p}\right]\right)\left\| \textbf{K}\left[1-e^{-h}\right]\right\|^{n-1},
\end{align*}
the last inequality taking place according to \eqref{ineg norme}. We finally have the following inequality
\begin{align*}
	&\left|\Tr\Big( \textbf{K}\left[1-e^{-h}\right]^n\Big)-\Tr\Big(\textbf{K}\left[1-e^{-h_p}\right]^n\Big)\right| \\
	&\hspace{5.5cm}\leq  \Tr \left(\textbf{K}\left[1-e^{-h}\right]-\textbf{K}\left[1-e^{-h_p}\right]\right) n \left\| \textbf{K}\left[1-e^{-h}\right]\right\|^{n-1}\\
	&\hspace{5.5cm}\underset{p\rightarrow +\infty}{\longrightarrow}0.
\end{align*}
Thus,
\begin{align*}
\underset{p\rightarrow +\infty}{\lim}\exp\left(-\sum_{n= 1}^{\infty} \frac{1}{n}\Tr\Big(\mathbf{K}\left[1-e^{-h_p}\right]^n\Big)\right)&=\exp\left(-\underset{p\rightarrow +\infty}{\lim}\sum_{n= 1}^{\infty} \frac{1}{n}\Tr\Big(\mathbf{K}\left[1-e^{-h_p}\right]^n\Big)\right)\\
&=\exp\left(-\sum_{n= 1}^{\infty} \frac{1}{n}\underset{p\rightarrow +\infty}{\lim}\Tr\Big(\mathbf{K}\left[1-e^{-h_p}\right]^n\Big)\right)\\
&=\exp\left(-\sum_{n= 1}^{\infty} \frac{1}{n}\Tr\Big(\mathbf{K}\left[1-e^{-h}\right]^n\Big)\right).
\end{align*}\CQFD

\end{Proof}
\vspace{1cm}

In order to prove the convergence in distribution of $n(\rho)^{-1}\widetilde{M}_\rho(\mu)$, for $\mu\in{\cal M}_{\alpha,\beta}^+$, we study the convergence of its Laplace transform : for $\theta\geq 0$,
\begin{align*}
\ee\Big[\exp&\big(-\theta n(\rho)^{-1}\widetilde{M}_\rho(\mu)\big)\Big]
=\exp\big(\theta\ee[n(\rho)^{-1}M_\rho(\mu)]\big)\
\ee\Big[\exp\big(-\theta n(\rho)^{-1}M_\rho(\mu)\big)\Big]\\
&\hspace{-0.5cm}=\exp\big(\theta\ee[n(\rho)^{-1}M_\rho(\mu)]\big)\
\ee\left[\exp\left(-\int_{\real^d\times\real_+\times\real_+}\theta n(\rho)^{-1}m\mu(B(x,r))\Phi_\rho(dx,dr,dm)\right)\right].
\end{align*}
To compute this last term, we use Proposition \ref{Laplace}. The hypothesis \eqref{cond:integrability} in this proposition is satisfied because in our context of weighted balls model, we have $h(x,r,m)=m\mu(B(x,r))$ and therefore: 
\begin{align*}
	\int_{\real^d\times\real_+\times\real_+}(1-&e^{- m\mu(B(x,r))})K_\rho(x,x)f(r/\rho)g(m)dx\frac{dr}{\rho} dm \\
	&\leq \lambda (\rho)\underset{x\in\real^d}{\sup}K(x,x) \int_{\real^d\times\real_+\times\real_+} m\mu(B(x,r))f(r/\rho)g(m)dx\frac{dr}{\rho}dm \\
&\leq \lambda (\rho)\rho^d  v_d \mu(\real^d) \underset{x\in\real^d}{\sup}K(x,x) \left(\int_{\real_+}mg(m)dm\right)\left(\int_{\real_+}r^df(r)dr\right) <+\infty .
\end{align*}
Moreover, we have to check that $\disp \widehat{K}_\rho\left[1-e^{-h} \right]  \in L^2(\real^d\times\real_+\times\real_+)$, if we denote by $h$ the function given by $h(x,r,m)=m\mu(B(x,r))$ defined on $\real^d\times\real_+\times\real_+$.
\begin{align*}
&\int_{(\real^d\times\real_+\times\real_+)^2} \widehat{K}_\rho\left[1-e^{-h} \right]^2((x,r,m),(y,s,m'))dxdrdmdydsdm' \\
&= \int_{(\real^d\times\real_+\times\real_+)^2}\left(1-e^{- m\mu(B(x,r))}\right)g(m)\frac{f(r/\rho)}{\rho}K_\rho^2(x,y)\\
&\hspace*{7cm}\times\frac{f(s/\rho)}{\rho}g(m') \left(1-e^{- m'\mu(B(y,s))}\right)  dxdrdmdydsdm'\\
&\underset{(a)}{\leq} \lambda(\rho)^2\left(\underset{x\in\real^d}{\sup}K(x,x)\right)^2 \left(\int_{\real^d\times\real_+\times\real_+}\left(1-e^{- m\mu(B(x,r))}\right)g(m)\frac{f(r/\rho)}{\rho} dxdrdm\right)^2 \\
&\leq  \lambda(\rho)^2\left(\underset{x\in\real^d}{\sup}K(x,x)\right)^2 \left( \rho^d  v_d \mu(\real^d)\left(\int_{\real_+}mg(m)dm\right)\left(\int_{\real_+}r^df(r)dr\right) \right)^2 <+\infty,
\end{align*}
where inequality (a) stands because $\disp K_{\rho}^2(x,y)\leq K_{\rho}(x,x)K_{\rho}(y,y)\leq \lambda(\rho)^2\left(\underset{x\in\real^d}{\sup}K(x,x)\right)^2$ thanks to \eqref{eq:lambdarho2}.\\\\
We can now apply Proposition \ref{Laplace} to obtain:
\begin{eqnarray} 
\nonumber
\displaystyle \ee \left[\exp\left(-\int_{\real^d\times\real_+\times\real_+}\theta n(\rho)^{-1}m\mu(B(x,r))\Phi_\rho(dx,dr,dm)\right)\right]\\
\label{eq:laplace-field}
& \hspace{-4.5cm} =\displaystyle\exp\left(-\sum_{n= 1}^{\infty} \frac{1}{n}\Tr\Big(\mathbf{K}\left[1-e^{-\theta n(\rho)^{-1} h}\right]^n\Big)\right).
\end{eqnarray}

\noindent The Laplace transform of $n(\rho)^{-1}\widetilde{M}_\rho(\mu)$ thus rewrites
\begin{align}\label{eq:LT2lb}
&\ee\Big[e^{-\theta n(\rho)^{-1}\widetilde{M}_\rho(\mu)}\Big]=\exp\bigg(\int_{\real^d\times\real_+\times\real_+}\psi\big(\theta n(\rho)^{-1}m\mu(B(x,r))\big)\ \lambda(\rho) K(0)\frac{f(r/\rho)}{\rho}g(m) dxdrdm\bigg)\nonumber\\
&\hspace*{6cm}\times \exp\bigg(-\sum_{n\geq 2} \frac{1}{n}\Tr\Big( {\bf \widehat{K}_\rho}\Big[1-e^{-\theta n(\rho)^{-1}h}\Big]^n\Big) \bigg)
\end{align}
with $\psi (u)=e^{-u}-1+u$.\\
The convergence of \eqref{eq:LT2lb} derives from the following lemmas. The complete proofs of these lemmas are given in \cite{BCG2018}. 
\begin{lemme}
\label{lemme:trace2n}
For all $n\geq 2$, we have 
\begin{equation*}
\label{eq:Tr2n}
\mathrm{Tr}\Big({\bf \widehat{K}_\rho}\big[1-e^{-\theta n(\rho)^{-1}h}\big]^n\Big)
\leq \mathrm{Tr}\Big({\bf \widehat{K}_\rho}\big[1-e^{-\theta n(\rho)^{-1}h}\big]^2\Big)^{n/2}.
\end{equation*}
\end{lemme}
\begin{Proof}
The key point is that ${\bf \widehat{K}_\rho}\big[1-e^{-\theta n(\rho)^{-1}h}\big]$ is trace-class, and as a consequence a Hilbert-Schmidt operator thanks to Lemma \ref{trace-compact} so the proof of Lemma 2.9 in \cite{BCG2018} applies in the same way.
\CQFD
\end{Proof}
\begin{lemme}
\label{lemme:trace2}
Assume Conditions \eqref{eq:Kxx}, \eqref{eq:cfvol} and \eqref{eq:lambdarho2}, and consider $\mu \in{\cal M}_{\alpha,\beta}^+$.
Then there is a constant $M\in ]0,+\infty[$ such that,   
$$
\mathrm{Tr}\Big({\bf \widehat{K}_\rho}\big[1-e^{-\theta n(\rho)^{-1}h}\big]^2\Big)\leq
 M\theta^2\frac{\lambda(\rho)\rho^q}{n(\rho)^2}.  
$$
\end{lemme}
\begin{Proof}
The computations are analogous to that in \cite{BCG2018} for the model without weight. It is important to observe that the key point, namely inequality (25) in \cite{BCG2018}, remains true because $\mu\in\mathcal{M}_{\alpha,\beta}^+$ so thanks Proposition 2.2 (iii) in \cite{BD2009}, $\mu\in\mathcal{M}_{2,\beta}^+=\mathcal{M}_{\beta}^+$ using the notations of \cite{BCG2018}. To be complete, the constant $M$ is equal to $\disp C_KC_\mu C_f \left(\int_{\real_+}mg(m)dm\right)^2$ with the notations of Lemma 2.10 in \cite{BCG2018}.
\CQFD 
\end{Proof}

\noindent \begin{Proof} We give a short proof of Theorem \ref{main theorem}. In this non-stationary case, the proof follows the same general strategy as in \cite{BCG2018}. 
Roughly speaking, the limits are driven by the term $n=1$ in \eqref{eq:laplace-field} while the other terms ($n\geq 2$) are still negligible. \\
Note that, in this non-stationary setting,  the Poissonian limits for $n=1$ come from Theorem 1, Theorem 2 and Theorem 3 in \cite{Gobard2014} taking $f(x,r)=K(x,x)f(r)$ in our situation.\\ 
As in \cite{BCG2018}, it is now enough to show now that
$$\underset{\rho\rightarrow 0}{\lim} ~\Tr\Big({\bf \widehat{K}_\rho}\big[1-e^{-\theta n(\rho)^{-1}h}\big]^2\Big)=0$$
in the three regimes.\\
As a consequence of Lemma \ref{lemme:trace2}, it remains to show that $\frac{\lambda(\rho)\rho^q}{n(\rho)^2}\underset{\rho \rightarrow 0}{\longrightarrow}0$.
\begin{enumerate}[(i)]
\item \textit{Large-balls scaling}. Since $\underset{\rho \rightarrow 0}{\lim} ~\lambda(\rho)\rho^\beta=+\infty$, for $\rho$ small enough we have $\lambda(\rho)\rho^\beta\geq 1$ and so $\left(\lambda(\rho)\rho^\beta\right)^{1/\alpha}\geq \left(\lambda(\rho)\rho^\beta\right)^{1/2}$ with $\alpha\in ]1,2]$. Thus since $q>\beta$ we have :
\[0\leq \frac{\lambda(\rho)\rho^q}{n(\rho)^2}\leq\frac{\lambda(\rho)\rho^q}{\lambda(\rho)\rho^\beta}=\rho^{q-\beta}\underset{\rho \rightarrow 0}{\longrightarrow}0.\]
\item \textit{Intermediate scaling}. In this case, $n(\rho)=1$ and since $q>\beta$ we have :
\[0\leq \frac{\lambda(\rho)\rho^q}{n(\rho)^2}=\lambda(\rho)\rho^q=\lambda(\rho)\rho^\beta \rho^{q-\beta}\underset{\rho \rightarrow 0}{\longrightarrow}0.\]
\item \textit{Small-balls scaling}. Since we consider $\mu \in L^1(\real^d)\cap L^2(\real^d) $, Proposition 2.5-(ii) in \cite{BCG2018} (which remains correct for $\mu\in L^1(\real^d)\cap L^2(\real^d)$) ensures that we can take $q=2d$ and then with $n(\rho)=(\lambda(\rho)\rho^\beta)^{1/\gamma}$ since $\beta<\alpha d\leq 2d$ we have :
\[0\leq \frac{\lambda(\rho)\rho^{2d}}{ n(\rho)^2}
=\lambda(\rho)^{(\beta-2d)/\beta}\underset{\rho \rightarrow 0}{\longrightarrow}0.\]~\CQFD
\end{enumerate}
\end{Proof}


\section{Appendix: Lemmas for the Laplace transform of DPP} 
\label{sec:laplacelemme}

In this appendix, we state and prove three different lemmas used in section \ref{sec:proof} to prove Proposition \ref{Laplace}.

\begin{lemme}\label{ineg norme}
If $f, g$ are two real functions on $E$ such that $0\leq f \leq g$, then we have
\[\left\|\mathbf{K}\left[f \right]  \right\|\leq \left\|\mathbf{K}\left[g  \right] \right\|\]
where $\mathbf{K}\left[f \right]$ is the operator with kernel $K\left[f \right](x,y)=\sqrt{f(x)}K(x,y)\sqrt{f(y)}$ for $x,y\in E$.
\end{lemme}
\begin{Proof}
Recall that:
\[\left\|\textbf{K}\left[f \right] \right\|=\underset{h \in L^2(E)}{\sup}\frac{\langle \textbf{K}\left[f \right](h),h \rangle}{\left\|h\right\|_2^2}.\]
Let $h\in L^2(E)$.
\begin{align*}
\langle \textbf{K}\left[f \right](h),h \rangle &=\int_E \textbf{K}\left[f\right](h)(x)h(x)dx \\
&=\int_E\int_E \sqrt{f(x)}K(x,y)\sqrt{f(y)}h(y)h(x)dydx\\
&=\int_E\int_E \sqrt{\frac{f(x)}{g(x)}}\sqrt{g(x)}K(x,y)\sqrt{g(y)}\sqrt{\frac{f(y)}{g(y)}}h(y)h(x)dydx \\
&=\langle \textbf{K}\left[g \right]\left(l h\right), l h\rangle
\end{align*}
where $\displaystyle l=\sqrt{f/g}\leq 1$.
So we have:
\[\langle \textbf{K}\left[f \right](h),h \rangle \leq \left\|\textbf{K}\left[g \right] \right\| \left\|lh\right\|_2^2 \leq \left\|\textbf{K}\left[g\right] \right\| \left\|h\right\|_2^2\]
and the result follows.
\CQFD 
\end{Proof}

\noindent We recall the following result from Proposition 280 in \cite{Goss2011}:
\begin{lemme}\label{trace-compact}
A trace-class operator $\textbf{K}$ on a Hilbert space $\mathcal{H}$ is a compact operator.
\end{lemme}


\vspace{1cm}



\end{document}